\documentclass[preprint,12pt]{elsarticle}

\usepackage{amssymb,amsthm,amsmath,amsfonts,mhequ,color}

\journal{Stochastic Processes and their Applications}

\def\begmat{\left(\begin{array}}\def\endmat{\end{array}\right)}

\def\bi{\begin{itemize}\setlength{\itemsep}{0pt}} \def\ei{\end{itemize}}

\def\bl{\begin{list}{\labelitemi}{\leftmargin=1em}\setlength{\itemsep}{-2.5pt}}  \def\el{\end{list}}
\def\bn{\begin{enumerate}} \def\en{\end{enumerate}}
\def\bt{\begin{table}[h]} \def\et{\end{table}}
\def\bc{\begin{center}} \def\ec{\end{center}}
\def\T{{ \mathrm{\scriptscriptstyle T} }}

\newcommand{\bfpsi} {\mbox{\boldmath $\psi$}}

\newcommand{\bfx} {\mbox{\boldmath $x$}}

\newcommand{\abs}[1]{\left\vert#1\right\vert}

\newcommand{\norm}[1]{\left\Vert#1\right\Vert}

\newcommand \bbP{\mathbb{P}}
\newcommand \bbE{\mathbb{E}}

\newtheorem{theorem}{Theorem}[section]
\newtheorem{lemma}[theorem]{Lemma}

\newtheorem{corollary}[theorem]{Corollary}
\newtheorem{remark}[theorem]{Remark}

\theoremstyle{plain}

\theoremstyle{plain}

\theoremstyle{remark}

\theoremstyle{plain}

\newcommand \be{\begin{equs}}
\newcommand \ee{\end{equs}}

\newcommand \N{\mathrm{N}}

\begin{document}

\begin{frontmatter}



\title{Sub-optimality of some continuous shrinkage priors}


\author{Anirban Bhattacharya, David B. Dunson, Debdeep Pati, Natesh S. Pillai}

\address{Texas A\& M University, Duke University, Florida State University, Harvard University}

\begin{abstract}
Two-component mixture priors provide a traditional way to induce sparsity in high-dimensional Bayes models.   However, several aspects of such a prior, including computational complexities in high-dimensions, interpretation of exact zeros and non-sparse posterior summaries under standard loss functions, has motivated an amazing variety of continuous shrinkage priors, which can be expressed as global-local scale mixtures of Gaussians.  
Interestingly, we demonstrate that many commonly used shrinkage priors,  including the Bayesian Lasso, do not have adequate posterior concentration in high-dimensional settings.  

\end{abstract}


\end{frontmatter}


\section{Introduction}
With the recent flurry of activities in high-throughput data, taking advantage of  sparsity to perform statistical inference is a common theme in situations where the number of model parameters ($p$) increases with the sample-size ($n$). In such scenarios, 
penalization methods \citep{buhlmann2011statistics} can yield a point estimate very quickly. There is a rich theoretical literature justifying the optimality properties of such penalization approaches  \cite{zhao2007model,van2008high,zhang2008sparsity,meinshausen2009lasso,raskutti2011minimax,negahban2010unified}, with fast algorithms \cite{efron2004least} and compelling applied results leading to their routine use.\par
On the other hand, statistical theory for characterizing the uncertainty of model parameters using penalization methods in high dimensions has received comparatively less attention \cite{nickl2013confidence}. Bayesian approaches provide a natural measure of uncertainty through the induced posterior distribution. Most penalization methods have a Bayesian counterpart. For instance, $\ell_1$ and $\ell_2$ regularization methods are equivalent to placing zero-mean double-exponential and Gaussian priors respectively  on the parameter vector and the solutions of the corresponding optimization problems are precisely the mode of the Bayesian posterior distribution. Moreover, a Bayesian approach has distinct advantages in terms of choice of tuning parameters, allowing key penalty parameters to be marginalized over the posterior distribution instead of relying on cross-validation. Thus a fruitful line of research is to investigate the behavior of the entire posterior distribution of the Bayesian models corresponding to penalization methods.  \par
 The process of eliciting prior distributions can be very tricky in high-dimensions.  Two-component mixture priors with a point mass at zero are traditionally used in high-dimensional settings because of  their ability to produce exact zeros and ease of 
 eliciting hyperparameters based on the prior knowledge about the level of sparsity and the size of the signal coefficients. In  \cite{castilloneedles,castillo2014bayesian}, the authors showed optimality properties for carefully chosen point mass mixture priors in high-dimensional settings.  Recently,  in an insightful article \cite{polson2010shrink}, several arguments were raised against the point mass priors concerning interpretation of exact zeros and computational issues arising from exploring a very high-dimensional model space.  This prompted the authors of \cite{polson2010shrink} to seek for continuous analogues of point mass priors based on Gaussian scale mixtures which obviates the need to search over the huge model space.   These scale mixtures of Gaussian priors are designed to have a sharp peak near zero with heavy tails so as to emulate the point mass mixture priors. 
In the last few years, a huge variety of shrinkage priors have been proposed in the Bayesian literature \cite{park2008bayesian,tipping2001sparse,griffin2010inference,carvalho2010horseshoe,armagan2013generalized}. In \cite{armagan2013generalized} the authors studied shrinkage priors and provided simple sufficient conditions for posterior consistency in $p \le n$ settings.  However, results on quantifying posterior concentration using continuous shrinkage priors are scarce. 

  
Even from a purely practical point of view, considerable difficulties have arisen when attempts have been made to reflect prior beliefs on sparsity through the associated hyperparameters of these distributions.  For example, suppose we wish to estimate $\theta_0 \in \mathbb{R}^n$ from $y \sim \mbox{N}_n(\theta_0, I_n)$ under the prior knowledge that only a fraction of the coordinates of $\theta_0$ are non-zero.  What are the appropriate parameters one should choose in the Bayes Lasso formulation  \citep{park2008bayesian} to ensure efficient estimation of $\theta_0$? A first step towards answering such questions is to understand the concentration of shrinkage priors around sparse vectors. This is critically important in two aspects. First, optimal prior concentration is almost necessary for optimal posterior contraction rates under a variety of loss functions. Second, studying  the concentration of shrinkage priors around sparse vectors will yield insights into the geometry of shrinkage priors which can then be harnessed for prior elicitation for a broad class of models. \par

Our contribution in this paper is two fold. First, we obtain sharp bounds for the concentration of continuous shrinkage priors around sparse high dimensional vectors. This is quite challenging because the joint distributions of such priors  obtained through integrating several latent hyperparameters are often unwieldy to work with. One of the reasons why the point mass priors enjoy theoretical optimality properties is because they have optimal concentration around sparse vectors \citep{castilloneedles,castillo2014bayesian}. We show that the concentration of some of the commonly used continuous shrinkage priors can sometimes be smaller than that of the point mass priors by several orders of magnitude.
Second, using these results, we show that for the normal means problem, the Bayesian Lasso \cite{park2008bayesian} has sub-optimal posterior concentration around sparse vectors under the $\ell_2$ loss.   
Although not rigorously proven, we suspect the sub-optimality to be actually by a genuine power of  the sample size, making this a serious concern against routine use of such priors. The negative results about shrinkage priors obtained in this paper contribute further to our understanding of the `statistical efficiency \textit{vs.} computational efficiency' trade off  in high dimensional estimation.  
Technically, we rely on exploiting the conditional Gaussian formulation of continuous shrinkage priors and use results for small-ball probability of Gaussian distributions.    Although we focus exclusively on the normal means problem, the ideas  in this paper are applicable to other models including linear regression and to estimating high dimensional covariance matrices using Gaussian latent factor models \cite{pati2014posterior}. \par  

One of the key insights emerging from this paper is that to induce an appropriate dependence structure in the scale parameters to ensure that the concentration is sufficiently high.  
In a recent article \cite{bhattacharya2014dirichlet}, we use this insight to construct a prior distribution on high dimensional vectors called the `Dirichlet-Laplace' prior which enjoys optimal posterior concentration rates.  Also, while completing this paper, we became aware of a recent article \cite{van2014horseshoe} in which the authors establish optimal posterior concentration for horseshoe priors in the normal means problem. Their techniques rely on explicit expressions for the posterior means and variances in the normal means model.  However, the hyperparameters of the horseshoe priors are chosen using an empirical Bayes approach.  The optimality of the fully Bayesian version of the horseshoe prior as in \cite{carvalho2010horseshoe} remains an interesting open problem. 

Proofs of technical results are provided in Section 6. Proofs of auxiliary lemmata stated in Section 5 are in a supporting document.

\section{Preliminaries} \label{sec:prelim}

Given sequences $a_n, b_n$, we denote $a_n = O(b_n)$ if there exists a global constant $C$ such that $a_n \leq C b_n$ and $a_n = o(b_n)$ if $a_n/b_n \to 0$ as $n \to \infty$.  For a vector $ x \in \mathbb{R}^r$,  $\norm{x}_2$ denotes its Euclidean norm. We will use $\Delta^{r-1}$ to denote the $(r-1)$-dimensional simplex $\{ x = (x_1, \ldots, x_{r})^{\T} : x_j \geq 0, \sum_{j=1}^{r} x_j = 1\}$. Further, let $\Delta_0^{r-1}$ denote $\{ x = (x_1, \ldots, x_{r-1})^{\T} : x_j \geq 0, \sum_{j=1}^{r-1} x_j \leq 1\}$.

For a subset $S \subset \{1, \ldots, n\}$, let $|S|$ denote the cardinality of $S$ and define $\theta_S = (\theta_j : j \in S)$ for a vector $\theta \in \mathbb{R}^n$. 
Denote $\mbox{supp}(\theta)$ to be the \emph{support} of $\theta$, the subset of $\{1, \ldots, n\}$ corresponding to the non-zero entries of $\theta$. 
Let $\ell_0[q; n]$ denote the subset of  $\mathbb{R}^n$ consisting of $q$-sparse vectors $\theta$ with $|\mbox{supp}(\theta)| \leq q$:
\begin{align*}
\ell_0[q;n] = \{ \theta \in \mathbb{R}^n ~:~ \#(1 \leq j \leq n : \theta_j \neq 0) \leq q\}.
\end{align*}
Let $\mbox{DE}(\tau)$ denote a zero mean double-exponential or Laplace distribution with density $f(y) = (2 \tau)^{-1} e^{- \abs{y}/\tau}$ for $y \in \mathbb{R}$.

\section{Concentration properties of global-local priors} \label{sec:conc_prop}

\subsection{\bf{Motivation}}

For a high-dimensional vector $\theta \in \mathbb{R}^n$, a natural way to incorporate sparsity in a Bayesian framework is to use point mass mixture priors 
\begin{align}\label{eq:point_mass}
\theta_j \sim (1 - \pi) \delta_0 + \pi g_{\theta}, \quad j = 1, \ldots, n ,
\end{align}
where $\pi = \mbox{Pr}(\theta_j \neq 0)$, $\bbE\{ |\mbox{supp}(\theta)|  \mid \pi\} = n \pi$ is the prior guess on model size (sparsity level), and $g_{\theta}$ is an absolutely continuous density on $\mathbb{R}$.  A beta prior on $\pi$ leads to an automatic multiplicity adjustment \citep{scott2010bayes}. In   \cite{castilloneedles}, the authors established that prior (\ref{eq:point_mass}) with an appropriate beta prior on $\pi$ and suitable tail conditions on $g_{\theta}$ leads to a frequentist minimax optimal rate of posterior contraction in the normal means setting.  We shall revisit the normal means problem in Section \ref{subsec:norm_means}.

As mentioned before, the authors in \cite{polson2010shrink} noted certain unappealing aspects of \eqref{eq:point_mass} including computational complexities in high-dimensions, interpretation of exact zeros and non-sparse posterior summaries under common loss functions.   This has motivated a variety of continuous shrinkage priors \citep{park2008bayesian,carvalho2010horseshoe,griffin2010inference,hans2011elastic,armagan2013generalized}, resembling the two-component priors facilitating computation and interpretability. 
Almost all such shrinkage priors can be represented as global-local (GL) mixtures of Gaussians \cite{polson2010shrink}, 
\begin{align}\label{eq:lg}
\theta_j \sim \mbox{N}(0, \psi_j \tau), \quad \psi_j \sim f, \quad \tau \sim g,  
\end{align}
where $\tau$ controls global shrinkage towards the origin while the local scales $\{ \psi_j \}$ allow deviations in the degree of shrinkage.  If $g$ puts sufficient mass near zero and $f$ is appropriately chosen, the hope is that GL priors in (\ref{eq:lg}) can serve as a good enough proxy for (\ref{eq:point_mass}) through a continuous density concentrated near zero with heavy tails. 

The normal scale mixture representation in \eqref{eq:lg} allows for conjugate updating of $\theta$ and $\psi$ in a block facilitating computation in high dimensions. Moreover, a number of frequentist regularization procedures such as ridge regression, lasso, bridge and elastic net correspond to posterior modes under GL priors with appropriate choices of $f$ and $g$. For example, one obtains a double-exponential prior corresponding to the popular $\ell_1$ or lasso penalty if $f$ has an exponential distribution. However, unlike variable selection priors (\ref{eq:point_mass}), many aspects of shrinkage priors are poorly understood.  For example, even basic properties, such as the appropriate choices of $f$ and $g$ to have adequate prior concentration around a sparse vector are unknown till date. 


There has been a recent awareness of these issues, motivating a basic assessment of the marginal properties of shrinkage priors for a single $\theta_j$. Recent priors such as the horseshoe \cite{carvalho2010horseshoe} and generalized double Pareto \cite{armagan2013generalized} are carefully formulated to obtain marginals having a high concentration around zero with heavy tails.  This is well justified, but as we will see below, such marginal behavior alone is not sufficient; it is necessary to study the joint distribution of $\theta$ on $\mathbb{R}^n$.  Specifically, we recommend studying the prior concentration $\bbP(\norm{\theta - \theta_0} < t_n)$ where the true parameter $\theta_0$ is assumed to be sparse: $\theta_0 \in \ell_0[q_n; n]$ with the number of non-zero components $q_n \ll n$ and
\begin{align}\label{eqn:tn}
t_n = n^{\delta/2} \quad \mathrm{with} \quad \delta \in (0,1).  
\end{align}

In models where $q_n \ll n$,  the prior must place sufficient mass around sparse vectors to allow for good posterior contraction; see Section \ref{subsec:norm_means} for further details. Now, as a first illustration, consider the following two extreme scenarios: i.i.d. standard normal priors for the individual components $\theta_j$ \textit{vs.} point mass mixture priors given by \eqref{eq:point_mass}.


\begin{theorem}\label{thm:basic_stmt}
Assume that $\theta_0 \in \ell_0[q_n; n]$ with $q_n = o(n)$. Then, for i.i.d standard normal priors on $\theta_j$, 
\begin{align}\label{thm:norm_conc}
\bbP(\norm{\theta - \theta_0}_2 < t_n) \leq e^{-c\, n}. 
\end{align}
For point mass mixture priors \eqref{eq:point_mass} with $\pi \sim \mbox{Beta}(1, n + 1)$ and $g_{\theta}$ being a standard Laplace distribution $g_{\theta} \equiv \mathrm{DE}(1)$, 
\begin{align}\label{thm:pm_conc}
\bbP(\norm{\theta - \theta_0}_2 < t_n) \geq e^{-c \max\{q_n, \norm{\theta_0}_1\}}.
\end{align}
\end{theorem}
\begin{proof}
Using $\norm{\theta}_2^2 \sim \chi_n^2$, the claim made in \eqref{thm:norm_conc} follows from an application of a Gaussian small probability in Lemma \ref{lem:anderson} and standard chi-square deviation inequalities. In particular, the exponentially small concentration also holds for $\bbP(\norm{\theta_0}_2 < t_n)$. The second claim \eqref{thm:pm_conc} follows from results in \cite{castilloneedles}. 
\end{proof}
As seen from Theorem \ref{thm:basic_stmt}, the point mass mixture priors have much improved concentration around sparse vectors, as compared to the i.i.d. normal prior distributions. The theoretical properties enjoyed by the point mass mixture priors can partly be attributed to this improved concentration. The above comparison suggests that it is of merit to evaluate a shrinkage prior in high dimensional models under sparsity assumption by obtaining its concentration rates around sparse vectors. 
The main results are given in Section \ref{sec:main_pf}. Recall the GL priors presented in \eqref{eq:lg} and the sequence $t_n$ in \eqref{eqn:tn}.

\subsection{\bf{Prior concentration for global priors}} This simplified setting involves  only a global parameter, \textit{i.e.}, $\psi_j = 1$ for all $j$. This subclass includes the important example of ridge regression, with $\tau$ routinely assigned an inverse-gamma prior, $\tau \sim \mbox{IG}(\alpha, \beta)$.  
\begin{theorem}\label{thm:conc_invgam}
Assume $\theta \sim \mathrm{GL}$ with $\psi_j = 1$ for all $j$. If the prior $g$ on the global parameter $\tau$ has an $\mathrm{IG}(\alpha, \beta)$ distribution, then
\begin{align}\label{eq:conc_ig}
\bbP(\norm{\theta}_2 < t_n) \leq e^{- C n^{1- \delta} },
\end{align}
where $C > 0$ is a constant depending only on $\alpha$ and $\beta$ and $\delta$ is from \eqref{eqn:tn}.  
\end{theorem}
The above theorem shows that compared to i.i.d. normal priors \eqref{thm:norm_conc}, the prior concentration does not improve much under an inverse-gamma prior on the global variance regardless of the hyperparameters (provided they don't scale with $n$) even when $\theta_0=0$.  Concentration around $\theta_0$ away from zero will clearly be even worse.  Hence, such a prior is not well-suited in high-dimensional settings, confirming empirical observations  documented in \cite{gelman2006prior,polson2011half}. It is also immediate that the same concentration bound in \eqref{eq:conc_ig} would be obtained for the giG family of priors on $\tau$.

In \cite{polson2011half}, the authors instead recommended a half-Cauchy prior as a default choice for the global variance (also see \cite{gelman2006prior}).  We consider the following general class of densities on $(0, \infty)$ for $\tau$, to be denoted $\mathcal{G}$ henceforth, that satisfy:  (i) $g(\tau) \leq M$ for all $\tau \in (0, \infty)$  (ii) $g(\tau) > 1/M$ for all $\tau \in (0, 1)$, for some constant $M > 0$. Clearly, $\mathcal{G}$ contains the half-Cauchy and exponential families. The following result provides concentration bounds for these priors. 

\begin{theorem}\label{thm:conc_g}
Let $\norm{\theta_0}_2 = o(\sqrt{n})$. Recall $t_n$ and $\delta$ from \eqref{eqn:tn}. If the prior $g$ on the global parameter $\tau$ belongs to the class $\mathcal{G}$ above then,
\begin{align}
& C_1 e^{-(1-\delta) \log n} \leq \bbP(\norm{\theta}_2 <  t_n ) \leq C_2 e^{- (1- \delta) \log n} \label{eq:conc_gonly_c}.
\end{align}
Furthermore, if $\frac 14 \norm{\theta_0}_2 > {t_n}$, then
\begin{align}
e^{- c_1 n \log a_n }\leq \bbP(\norm{\theta - \theta_0}_2 <  t_n ) \leq  e^{- c_2 n \log a_n }, \label{eq:conc_gonly_nc}
\end{align}
where $a_n = \norm{\theta_0}_2/\sqrt{2} t_n > 1$ and $c_i, C_i > 0$ are constants with $C_1, C_2$ depending only on $M$ in the definition of $\mathcal{G}$ and $c_1$ depending on $M$ and $\delta$. 
\end{theorem}

Equation \eqref{eq:conc_gonly_c} in Theorem \ref{thm:conc_g} shows that the prior concentration around zero can be dramatically improved from exponential to polynomial with a careful prior on $\tau$ that can assign sufficient mass near zero, such as the half-Cauchy prior \cite{gelman2006prior, polson2011half}.  Unfortunately, as \eqref{eq:conc_gonly_nc} shows,  for signals of large magnitude one again obtains an exponentially decaying probability.  Hence, Theorem \ref{thm:conc_g} conclusively shows that global shrinkage priors are simply not flexible enough for high-dimensional problems. 

An inspection of the proof of both Theorems will reveal that the condition $t_n = n^{\delta/2}$ is only used at the last step to present the bound in its simplest form. A similar bound can be derived for other sequences $t_n$, for example, when $t_n$ grows logarithmically, a fact which is used in the proof of Theorem \ref{thm:lb_g} later. 

\subsection{\bf{Prior concentration for a class of GL priors}} Proving concentration results for the GL family \eqref{eq:lg} in the general setting presents a much harder challenge compared to Theorem \ref{thm:conc_g} since we now have to additionally integrate over the $n$ local parameters $\psi = (\psi_1, \ldots, \psi_n)$.  We focus on an important sub-class in Theorem \ref{thm:conc_lg} below, namely the exponential $\mbox{Exp}(\lambda)$ family for the distribution of $g$ in \eqref{eq:lg}. For analytical tractability, we additionally assume that $\theta_0$ has only one non-zero entry. The interest in the $\mbox{Exp}(\lambda)$ arises from the fact that normal-exponential scale mixtures give rise to the double-exponential family \cite{west1987scale}: $\theta \mid \psi \sim N(0, \psi \sigma^2), \psi \sim \mbox{Exp}(1/2)$ implies $\theta \sim \mbox{DE}(\sigma)$, and hence this family of priors can be considered as a Bayesian version of the lasso \cite{park2008bayesian}. We now state a concentration result for this class. The proof is provided in the supporting document.
\begin{theorem}\label{thm:conc_lg}
Assume $\theta \sim \mathrm{GL}$ with $g \in \mathcal{G}$ and $f \equiv \mathrm{Exp}(\lambda)$ for some constant $\lambda > 0$. Also assume $\theta_0$ has only one non-zero entry.  Then, for a global constant $C_1 > 0$ depending only on $M$ in the definition of $\mathcal{G}$, for all $\sqrt{w_n} = t_n < \norm{\theta_0}_2/4$,
 \begin{align}\label{eq:lg_ub_stmt}
& \bbP(\norm{ \theta - \theta_0}_2 \leq t_n) \leq C_1 \, \int_{\psi_1 = 0}^{\infty} \frac{\psi_1^{(n-3)/2}}{ \big\{ \psi_1 + \norm{\theta_0}_2^2/(\pi w_n) \big \}^{(n-3)/2}  } e^{-\psi_1} d \psi_1 . 
\end{align}
Let $v_n = r_n^2$ satisfy $v_n = O(\sqrt{n})$. Then, for $\norm{\theta_0}_2 \geq 1/\sqrt{n}$, 
\begin{align}\label{eq:lg_lb_stmt}
& \bbP(\norm{ \theta - \theta_0}_2 \leq r_n) \geq C_2 e^{-d_2 \sqrt{n}} \, \int_{\psi_1 =c_1 \norm{\theta_0}_2^2}^{\infty} \frac{\psi_1^{(n-3)/2}}{ \big\{ \psi_1 + \norm{\theta_0}_2^2/(\pi v_n) \big \}^{(n-3)/2}  } e^{-\psi_1} d \psi_1,
\end{align}  
where $c_1, d_2, C_2$ are positive global constants with $c_1 \geq 2$ and $C_2$ depends only on $M$ in the definition of $\mathcal{G}$. 
\end{theorem}
A more interpretable corollary can be stated as follows. 
\begin{corollary}\label{cor:thm_lg}
Assume $\theta \sim \mathrm{GL}$ with $g \in \mathcal{G}$ and $f \equiv \mathrm{Exp}(\lambda)$ for some constant $\lambda > 0$. Also assume $\theta_0$ has only one non-zero entry and $\norm{\theta_0}_2^2 > \log n$. Then, for a global constant $C > 0$ depending only on $M$ in the definition of $\mathcal{G}$,
\begin{align}\label{thm:lg_ub_cor}
\bbP(\norm{\theta - \theta_0}_2 < \sqrt{\log n}) \leq e^{- C \sqrt{n}} . 
\end{align}
\end{corollary}
Corollary \ref{cor:thm_lg} asserts that even in the simplest deviation from the null model with only one signal, one continues to have exponentially small concentration under an exponential prior on the local scales. From \eqref{thm:pm_conc} in Theorem \ref{thm:basic_stmt}, appropriate point mass mixture priors \eqref{eq:point_mass} would have $\bbP( \norm{\theta - \theta_0}_2 < t_n) \geq e^{ - C \norm{\theta_0}_1}$ under the same conditions as above, clearly showing that the wide difference in concentration still persists.


\section{\bf{Posterior lower bounds in normal means}}\label{subsec:norm_means}

We have discussed the prior concentration for a high-dimensional vector $\theta$ without alluding to any specific model so far.
In this section we show how prior concentration impacts posterior inference for  the widely studied normal means problem\footnote{Although we study the normal means problem, the ideas and results in this section are applicable to other models such as non-parametric regression and factor models.} (see \cite{donoho1992maximum,johnstone2004needles,castilloneedles} and references therein):  
\begin{align}\label{eq:norm_means}
 y_i &= \theta_i + \epsilon_i, \quad \epsilon_i \sim \mbox{N}(0, 1), \quad 1 \leq i \leq n. 
\end{align}

The minimax rate $s_n$ for the above model is given by $s^2_n  \asymp q_n \log(n/q_n)$ when $\theta_0 \in \ell_0[q_n;n]$.
For this model \cite{castilloneedles} recently established that  for point mass priors for $\theta$ with $\pi \sim \mbox{beta}(1, \kappa n + 1)$ and $g_{\theta}$ having Laplace like or heavier tails,  the posterior contracts at the minimax rate, \textit{i.e.,}
$\bbE_{n, \theta_0} \bbP(\norm{\theta - \theta_0}_2 < M s_n \mid y) \to 1$ for some constant $M > 0$. Thus we see that carefully chosen point mass priors are indeed optimal\footnote{It is important that the hyper parameter for $\pi$ depends on $n$. We do not know if the result holds without this}. 

However not all choices for $g_\theta$ lead to optimal proceedures; \cite{castilloneedles} also showed that if $g_{\theta}$ is instead chosen to be standard Gaussian, \textit{the posterior does not contract at the minimax rate}, \textit{i.e.}, one could have $\bbE_{n, \theta_0} \bbP (\norm{\theta - \theta_0}_2 < s_n \mid y) \to 0$ for signals of sufficiently large magnitude. 
To establish such a posterior lower-bound result, building on the work of \cite{castillo2008lower}, \cite{castilloneedles} showed that given a fixed sequence $t_n$, if there exists a sequence $r_n$ ($r_n > t_n$) such that 
\begin{align}
\frac{\bbP(\norm{\theta - \theta_0}_2 < t_n )}{\bbP(\norm{\theta - \theta_0}_2 < r_n )} = o(e^{-r_n^2}), \label{eq:lb_ratio}
\end{align}
then $\bbP(\norm{\theta - \theta_0}_2 < t_n \mid y) \to 0$. This immediately shows the importance of studying the prior concentration. Intuitively, \eqref{eq:lb_ratio} would be satisfied when the prior mass of the bigger ball $\norm{\theta - \theta_0}_2 < r_n$ is almost entirely contained in the annulus with inner radius $t_n$ and outer radius $r_n$, so that the smaller ball $\norm{\theta - \theta_0}_2 < t_n$ barely has any prior mass compared to the bigger ball. As an illustrative example, in the i.i.d. $\mbox{N}(0, 1)$ example with $t_n = s_n$, setting $r_n = \sqrt{n}$ would satisfy \eqref{eq:lb_ratio} above, proving that i.i.d. $\mbox{N}(0, 1)$ priors are sub-optimal.
Our goal is to investigate whether a similar phenomenon persists for global-local priors in light of the concentration bounds developed in Theorems \ref{thm:conc_g} and \ref{thm:conc_lg}. 

As in Section 3.2, we first state our posterior lower bound result for the case where there is only a global parameter. 
\begin{theorem}\label{thm:lb_g}
Suppose we observe $y \sim \mbox{N}_n(\theta_0, \mathrm{I}_n)$ and \eqref{eq:norm_means} is fitted with a $\mathrm{GL}$ prior on $\theta$ such that $\psi_j = 1$ for all $j$ and the prior $g$ on the global parameter $\tau$ lies in $\mathcal{G}$.  Assume $\theta_0 \in \ell_0[q_n; n]$ where $q_n/n \to 0$ and $\norm{\theta_0}_2 > s_n$, with $s_n^2 = q_n \log(n / q_n)$ being the minimax squared error loss over $\ell_0[q_n; n]$. Then, 
$\bbE_{n, \theta_0} \bbP(\norm{\theta - \theta_0}_2 \leq \sqrt{A}s_n \mid y) \to 0$ for any constant $A > 0$.
\end{theorem}
\begin{proof}
Without loss of generality, assume $\norm{\theta_0}_2 = o(\sqrt{n})$, since the posterior mass with a prior centered at the origin would be smaller otherwise. Choosing $t_n = \sqrt{A}s_n$, $r_n$ to be a sequence such that $t_n < r_n < \norm{\theta_0}_2$ and resorting to the two-sided bounds in Theorem \ref{thm:conc_g}, the ratio in \eqref{eq:lb_ratio} is smaller than $(t_n/r_n)^n$, and hence $e^{r_n^2} (t_n/r_n)^n \to 0$ since $r_n \leq \norm{\theta_0}_2 = o(\sqrt{n})$. 
\end{proof}
 Theorem \ref{thm:lb_g} states that a GL prior with only a global scale is sub-optimal if $\norm{\theta_0}_2 > s_n$. Observe that in the complementary region $\{ \norm{\theta_0}_2 \leq s_n\}$, the estimator $\hat{\theta} \equiv 0$ attains squared error in the order of $q_n \log(n/q_n)$, implying the condition $\norm{\theta_0}_2 > s_n$ is hardly stringent. 
Clearly, the absence of local scales makes it challenging to estimate  both coefficients of different signal strengths simultaneously.

Next, we state a result for the sub-class of GL priors as in Theorem \ref{thm:conc_lg}, i.e., when $f$ has an exponential distribution leading to a double-exponential distribution marginally.  
\begin{theorem}\label{thm:lb_lg}
Suppose we observe $y \sim \mbox{N}_n(\theta_0, \mathrm{I}_n)$ and the model in \eqref{eq:norm_means} is fitted with a $\mathrm{GL}$ prior on $\theta$ such that $g$ lies in $\mathcal{G}$ and $f \equiv \mathrm{Exp}(\lambda)$ for some constant $\lambda > 0$. Assume $\theta_0 \in \ell_0[q_n; n]$ with $q_n = 1$ and $\norm{\theta_0}_2^2 /\log n \to \infty$. Then, $\bbE_{n, \theta_0} \bbP(\norm{\theta - \theta_0}_2 \leq \sqrt{A\log n} \mid y) \to 0$ for any constant $A > 0$. 
\end{theorem}
From \cite{castilloneedles}, appropriate point mass mixture priors would assign increasing mass with $n$ to the neighborhood $\{\norm{\theta - \theta_0}_2 \leq \sqrt{A\log n}\}$ for sufficiently large $A > 0$. Indeed, \cite{donoho1992maximum} established that the minimax rate for $\ell[1; n]$  is $\sqrt{2\log n}\{1+o(1)\}$. 
Hence, the Bayesian lasso \cite{park2008bayesian} is sub-optimal even in the simplest deviation from the null model with only one moderately sized signal. We believe that the conclusions would continue to be valid if one only assumes $f$ to have exponential tails plus some mild conditions on the behavior near zero. However, the assumptions of Theorem \ref{thm:lb_lg} precludes the case when $f$ has polynomial tails, such as the horseshoe \cite{carvalho2010horseshoe} and generalized double Pareto \cite{armagan2013generalized}. Very recently, \cite{van2014horseshoe} showed optimal posterior concentration for the horseshoe prior where the global parameter $\tau$ is estimated using an empirical Bayes approach. 

Another important question beyond the scope of the current paper should concern the behavior of the posterior when one plugs in an empirical Bayes estimator of the global parameter $\tau$. However, we show below that the ``optimal'' sample-size dependent plug-in choice $\tau_n = c^2/\log n$ (so that marginally $\theta_j \sim \mbox{DE}(c/\sqrt{\log n})$ ) for the lasso estimator \cite{negahban2010unified} produces a sub-optimal posterior:
\begin{theorem}\label{thm:lb_plugin}
Suppose we observe $y \sim \mbox{N}_n(\theta_0, \mathrm{I}_n)$ and \eqref{eq:norm_means} is fitted with a $\mathrm{GL}$ prior on $\theta$ such that $\tau$ is   deterministically chosen to be $\tau_n$, i.e., $g \equiv \delta_{\tau_n}$ for a non-random sequence $\tau_n$ and $f \equiv \mathrm{Exp}(\lambda)$ for some constant $\lambda > 0$.  Assume $\theta_0 \in \ell_0[q_n; n]$ with $q_n = 1$, $\norm{\theta_{0}} > 4s_n$, $\norm{\theta_0}\log n  = o(n)$ and $\tau_n = c/\log n$ is used as the plug-in choice. Then, $\bbE_{n, \theta_0} \bbP(\norm{\theta - \theta_0}_2 \leq s_n \mid y) \to 0$, with $s_n^2 = q_n \log(n / q_n)$ being the minimax squared error loss over $\ell_0[q_n; n]$. 
\end{theorem}

\section{Auxiliary results}\label{sec:aux}
We collect a set of Lemmata we need to prove the main results. All proofs are deferred to the supporting document.

An important tool used throughout is a two-sided Gaussian small ball probability \cite{van2008reproducing}.  
\begin{lemma}\label{lem:anderson}
Suppose $\theta \sim {\N}_n(0, \Sigma)$ with $\Sigma$ be a positive definite matrix and $\theta_0 \in \mathbb{R}^n$. Let $\norm{\theta_0}_{\mathbb{H}}^2 = \theta_0^{\T} \Sigma^{-1} \theta_0$. Then, for any $t > 0$,
\be \label{eqn:Andlb}
e^{- \frac{1}{2} \; \norm{\theta_0}_{\mathbb{H}}^2}  \bbP(\norm{\theta}_2 \leq t) \leq \bbP(\norm{\theta - \theta_0}_2 \leq t).
\ee
In addition, if $\Sigma$ is diagonal and $\norm{\theta_0}_0 = 1$, i.e., $\theta_0$ has only one non-zero entry, then  for all $t < \frac14 \norm{\theta_0}_2$, 
\be \label{eq:anderson_UB}
\bbP(\norm{\theta - \theta_0}_2 \leq t) \leq e^{- \frac{\norm{\theta_0}_{\mathbb{H}}^2}{4} } \bbP(\norm{\theta}_2 \leq t). 
\ee
If the condition $\norm{\theta_0}_0 = 1$ above is replaced by $\Sigma = \sigma^2 I$ for some $\sigma > 0$, then for all $t < \frac14 \norm{\theta_0}_2$, the bound in \eqref{eq:anderson_UB}  holds as well. \end{lemma}
\begin{remark}It is well known that among balls of fixed radius, a zero mean multivariate normal distribution places the maximum mass on the ball centered at the origin. Lemma \ref{lem:anderson} provides a sharp bound on the probability of shifted balls in terms of the centered probability and the size of the shift, measured via the RKHS norm $\norm{\theta_0}_{\mathbb{H}}^2$. 
\end{remark}

We next state Lemma \ref{lem:incomplete_gamma} to bound an incomplete gamma integral from below. Recall $\int_{\tau = 0}^{\infty} \tau^{-n/2} e^{-a_n/(2 \tau)} d\tau = \Gamma(n/2-1) (2/a_n)^{n/2-1}$. Lemma \ref{lem:incomplete_gamma} shows that the same integral over $(0, 1)$ is of the same order when $a_n \precsim n$.  

\begin{lemma}\label{lem:incomplete_gamma}
For a sequence $a_n \leq n/(2e)$, $\int_{\tau = 0}^1 \tau^{-n/2} e^{-a_n/(2\tau)} d\tau \geq (2/a_n)^{n/2-1} \Gamma(n/2-1) \xi_n$, where $\xi_n \uparrow 1$ with $(1 - \xi_n) \leq D/\sqrt{n}$ for some constant $D > 0$. 
\end{lemma}

We state a two sided bound for the complementary error function. 
\begin{lemma}\label{lem:erf}
Let $\mbox{erfc}(x) = \frac{2}{\sqrt{\pi}}\int_{x}^{\infty} e^{-t^2} dt$ denote the complementary error function. Then,
\begin{align}
& \sqrt{\pi} e^x \mathrm{erfc}(\sqrt{x}) \leq \frac{1}{\sqrt{x + 1/\pi}} \label{eq:erfc_ub} \\
& \sqrt{\pi} e^x \mathrm{erfc}(\sqrt{x}) \geq \bigg\{ \frac{1}{\sqrt{x}}  \bigg\}^{1 + \delta} \label{eq:erfc_lb}
\end{align}
where \eqref{eq:erfc_lb} holds for any $\delta > 0$ provided $x \geq 2$. 
\end{lemma}

We next state the Dirichlet integral formula (4.635 in \cite{gradshteyn1980corrected}) to simplify a class of integrals over the simplex $\Delta^{n-1}$:
\begin{lemma}\label{lem:dir_formula}
Let $h(\cdot)$ be a Lebesgue integrable function and $\alpha_j > 0, j = 1, \ldots, n$. Then,
\be
\int_{\sum x_j \leq 1} h\big(\sum x_j\big) \prod_{j=1}^n x_j^{\alpha_j - 1} d \bfx = \frac{\prod_{j=1}^n \Gamma(\alpha_j)}{\Gamma \big(\sum_{j=1}^n \alpha_j\big)} \int_{t=0}^1 h(t) \, t^{( \sum \alpha_j) - 1} dt. 
\ee
\end{lemma}
Lemma \ref{lem:dir_formula} follows simply by noting that the left hand side is $\bbE h(\sum_{j=1}^n X_j) $ up to normalizing constants where $(X_1, \ldots, X_n) \sim \mbox{Diri}(\alpha_1, \ldots, \alpha_n, 1)$, so that $\sum_{j=1}^n X_j \sim \mbox{Beta}(\sum \alpha_j, 1)$. Such probabilistic intuitions help us to reduce a more complicated integral stated in Lemma \ref{lem:dickey_integral} below. The proof of Lemma \ref{lem:dickey_integral} utilizes a beautiful identity found in \cite{dickey1968three}. We didn't find any reference for Lemma \ref{lem:dickey_integral}, though a related integral with $n/2$ in the exponent in the denominator appears in \cite{gradshteyn1980corrected}. 
\begin{lemma}\label{lem:dickey_integral}
Let $q_j, j = 0, 1, \ldots, n$ be positive numbers. Then, 
\begin{align*}
\int_{\sum x_j \leq 1 } \frac{\prod_{j=1}^n x_j^{-1/2}}{  [ \sum_{j=1}^n q_j x_j + q_0 ]^{n/2-1} } d \bfx = \frac{\Gamma(1/2)^n}{\Gamma(n/2)} q_0(n/2-1) \int_{x=0}^1 \frac{x^{n/2-2} (1-x)}{\prod_{j=1}^n \sqrt{ q_j x + q_0 }} dx .
\end{align*}
\end{lemma}


\section{Proofs of results in Sections 3 \& 4}\label{sec:main_pf}

In this section, we prove the main results of the paper other than Theorem \ref{thm:conc_lg}, which is deferred to the supporting document due to its length.

For GL shrinkage priors of the form (\ref{eq:lg}), given $\bfpsi = (\psi_1, \ldots, \psi_n)^{\T}$ and $\tau$, the elements of $\theta$ are conditionally independent with $\theta \mid \bfpsi, \tau \sim \mbox{N}_n(0, \Sigma)$ with $\Sigma = \mathrm{diag}(\psi_1 \tau, \ldots, \psi_n \tau)$. Hence we can use Lemma \ref{lem:anderson} to obtain
\be \label{eq:main_bd1}
& e^{- 1/(2 \tau) \sum_{j=1}^n \theta_{0j}^2/ \psi_j }  \, \bbP(\norm{\theta}_2 <  t_n \mid \bfpsi, \tau) \leq \bbP( \norm{\theta - \theta_0}_2 \leq t_n \mid \bfpsi, \tau).
\ee
Furthemore, if $\psi_j = 1$ for all $j$, then again by using Lemma \ref{lem:anderson} we obtain
\be
 \label{eq:main_bd1u}
  \bbP( \norm{\theta - \theta_0}_2 \leq t_n \mid \bfpsi, \tau)&\leq e^{- 1/(4 \tau) \sum_{j=1}^n \theta_{0j}^2 }\bbP(\norm{\theta}_2 <  t_n \mid \bfpsi, \tau). 
\ee
Letting $X_j = \theta_j^2$, $X_j$'s are conditionally independent given $(\bfpsi, \tau)$ with $X_j$ having a density $f(x_j \mid \bfpsi, \tau) = D /(\sqrt{\tau \psi_j x_j}) e^{-x_j/(2\tau \psi_j) }$ on $(0, \infty)$, where $D = 1/(\sqrt{2 \pi})$. Hence, with $w_n = t_n^2$, 
\begin{align}\label{eq:main_bd2}
\bbP(\norm{\theta}_2 < t_n \mid \bfpsi, \tau) =  D^n \int_{\sum x_j \leq w_n} \prod_{j=1}^n \frac{1}{\sqrt{x_j \tau \psi_j }} e^{-x_j/(2 \tau \psi_j )} d \bfx . 
\end{align}
For sake of brevity, we use $\{ \sum x_j \leq w_n \}$ in \eqref{eq:main_bd2} and all future references to denote the region $\{ \bfx \in \mathbb{R}^n : x_j \geq 0 \, \forall \, j = 1, \ldots, n, \, \sum_{j=1}^n x_j \leq w_n \}$. 
To estimate two-sided bounds for the marginal concentration $\bbP( \norm{\theta - \theta_0}_2 \leq t_n)$, we need to combine \eqref{eq:main_bd1} \& \eqref{eq:main_bd2} and integrate out $\bfpsi$ and $\tau$ carefully. We start by proving Theorem \ref{thm:conc_invgam} \& Theorem \ref{thm:conc_g} where one only needs to integrate out $\tau$. 

\subsection*{\bf{Proof of Theorem \ref{thm:conc_invgam}}}

In \eqref{eq:main_bd2}, set $\psi_j = 1$ for all $j$, recall $D  = 1/\sqrt{2 \pi}$ and $w_n = t_n^2$, and integrate over $\tau$ to obtain,
\begin{align}\label{eq:glob_marg_c}
\bbP(\norm{\theta}_2 \leq t_n ) 
= D^n \int_{\tau = 0}^{\infty} g(\tau) \bigg[  \int_{\sum x_j \leq w_n} \prod_{j=1}^n \frac{1}{\sqrt{x_j \tau}} e^{-x_j/(2 \tau)} d \bfx \bigg] d\tau .
\end{align} 
Substituting $g(\tau) = c \tau^{-(1+\alpha)} e^{-\beta/\tau}$ with $c = \beta^{\alpha}/\Gamma(\alpha)$ and using Fubini's theorem to interchange the order of integration between $x$ and $\tau$, \eqref{eq:glob_marg_c} equals
\begin{align}
& c D^n \int_{\sum x_j \leq w_n} \prod_{j = 1}^n \frac{1}{\sqrt{x_j}} \bigg [ \int_{\tau = 0}^{\infty} \tau^{-(1 + n/2+ \alpha) } e^{- \frac{1}{2 \tau}(2 \beta + \sum x_j)} d\tau \bigg] d \bfx \notag \\
& = c D^n 2^{n/2+\alpha} \Gamma(n/2+\alpha) \int_{\sum x_j \leq w_n} \frac{1}{(2 \beta + \sum x_j)^{n/2 + \alpha} } \prod_{j = 1}^n \frac{1}{\sqrt{x_j}}  \, d \bfx \notag \\
& = c D^n 2^{n/2+\alpha} w_n^{n/2} \Gamma(n/2+\alpha)  \int_{\sum x_j \leq 1} \frac{1}{(2 \beta + w_n \sum x_j)^{n/2 + \alpha} } \prod_{j = 1}^n \frac{1}{\sqrt{x_j}}  \, d \bfx. \label{eq:glob_ig_1}
\end{align}
Lemma \ref{lem:dir_formula} with $h(t) = 1/(2 \beta + w_n t)^{n/2 + \alpha}$ applied to \eqref{eq:glob_ig_1} implies
\begin{align}
\bbP(\norm{\theta}_2 \leq t_n )  & = c D^n 2^{n/2+\alpha}  w_n^{n/2} \Gamma(n/2+\alpha)  \frac{\Gamma(1/2)^n}{ \Gamma(n/2)} \int_{t=0}^1 \frac{t^{n/2-1}}{(2 \beta + w_n t)^{n/2 + \alpha}} dt.   \label{eq:glob_ig_2}
\end{align}
Substituting $D = 1/\sqrt{2 \pi}$, bounding $(2 \beta + w_n t)^{n/2 + \alpha} \geq (2 \beta)^{\alpha+1} (2 \beta + w_n t)^{n/2-1}$, and letting $\tilde{w}_n = w_n/(2 \beta)$, \eqref{eq:glob_ig_2} can be bounded above by
\begin{align*}
\frac{\Gamma(n/2+\alpha)}{\Gamma(n/2) \Gamma(\alpha) (2 \beta)^{\alpha +1 }} \tilde{w}_n^{n/2} \int_{t=0}^1 \frac{t^{n/2-1}}{(1 + \tilde{w}_n t)^{n/2-1}} dt  \leq \frac{w_n \Gamma(n/2+\alpha)}{\Gamma(n/2) \Gamma(\alpha) (2 \beta)^{\alpha +1 }}  \bigg(\frac{\tilde{w}_n}{ 1 + \tilde{w}_n} \bigg)^{n/2-1},
\end{align*}
where the second inequality above uses $t/(a+t)$ is an increasing function in $t > 0$ for fixed $a > 0$. By definition, $w_n = n^{\delta}$ for $0 < \delta < 1$ and hence $\frac{w_n \Gamma(n/2+\alpha)}{\Gamma(n/2) \Gamma(\alpha) (2 \beta)^{\alpha +1 }}$ can be bounded above by $e^{C_1 \log n}$. Also, using $(1-x)^{1/x} \leq e$ for all $x > 0$, $ \{ \tilde{w}_n/(1 + \tilde{w}_n) \}^{n/2-1}$ can be bound above by $e^{-  C_2 n/w_n} = e^{- C_2 n^{1 - \delta}}$. Hence the overall bound is $e^{- C n^{1 - \delta}}$ for some appropriate constant $C > 0$.  This completes the proof of Theorem \ref{thm:conc_invgam}. \qed 

\subsection*{\bf{Proof of Theorem \ref{thm:conc_g}}}

We start with the upper bound in \eqref{eq:conc_gonly_c}. The steps are similar as in the proof of Theorem \ref{thm:conc_invgam} and hence only a sketch is provided. 
Set $w_n = t_n^2$. Bounding $g(\tau) \leq M$ and interchanging order of integrals in \eqref{eq:glob_marg_c}, 
\begin{align}\label{eq:glob_gen_c1}
\bbP(\norm{\theta}_2 \leq t_n) \leq M D^n 2^{n/2-1} \Gamma(n/2-1) w_n \int_{\sum x_j \leq 1} \frac{1}{(\sum x_j)^{n/2-1} } \prod_{j = 1}^n \frac{1}{\sqrt{x_j}}  \, d \bfx .
\end{align}
Invoking Lemma \ref{lem:dir_formula} with $h(t) = (1/t)^{n/2 -1}$ in \eqref{eq:glob_gen_c1}, we obtain 
\small
\be 
M D^n 2^{n/2-1} \Gamma(n/2-1) w_n \frac{\Gamma(1/2)^n}{\Gamma(n/2)} \int_{x=0}^1 x^{n/2-1}/x^{n/2-1} dx = (M/2) \frac{w_n}{n/2 - 1} = C_2 n^{-(1-\delta)}.\\\label{eqn:globu}
\ee
\normalsize
From \eqref{eq:glob_gen_c1} and \eqref{eqn:globu} we obtain
\begin{align}\label{eq:glob_gen_c1}
\bbP(\norm{\theta}_2 \leq t_n) \leq C_2 n^{-(1-\delta)}
\end{align}
and thus the upper bound in \eqref{eq:conc_gonly_c} is proved. \par
We now turn towards proving the lower bound to the centered concentration in \eqref{eq:conc_gonly_c}. Recalling that $g(\tau) \geq 1/M$ on $(0, 1)$ for $g \in \mathcal{G}$, and interchanging integrals in \eqref{eq:glob_marg_c}, we have, with $K = 1/M$,
\begin{align}\label{eq:glob_c_lb1}
\bbP(\norm{\theta}_2 \leq t_n) \geq K D^n \int_{\sum x_j \leq w_n} \prod_{j = 1}^n \frac{1}{\sqrt{x_j}} \bigg [ \int_{\tau = 0}^{1} \tau^{-n/2 } e^{- \sum x_j / (2 \tau)} d\tau \bigg] d \bfx.
\end{align}
Clearly $\sum x_j \leq w_n$ and hence we can apply Lemma \ref{lem:incomplete_gamma} in \eqref{eq:glob_c_lb1} to get
\begin{align}\label{eq:glob_c_lb2}
\bbP(\norm{\theta}_2 \leq t_n) \geq K \xi_n D^n 2^{n/2-1} \Gamma(n/2-1) w_n \int_{\sum x_j \leq 1} \frac{1}{(\sum x_j)^{n/2-1} } \prod_{j = 1}^n \frac{1}{\sqrt{x_j}}  \, d \bfx .
\end{align}
The rest of the proof proceeds for showing the lower bound in \eqref{eq:conc_gonly_c} follows exactly as in the upper bound case from \eqref{eq:glob_gen_c1} onwards. 

Finally, we combine \eqref{eq:main_bd1} with \eqref{eq:main_bd1u} (with $\psi_j = 1$ for all $j$ in this case) to bound the non-centered probability in \eqref{eq:conc_gonly_nc}. For the upper bound, we additionally use $g(\tau) \leq M$ for all $\tau$ to obtain
\begin{align}
& \bbP( \norm{\theta - \theta_0}_2 \leq t_n) \leq  M D^n \int_{\sum x_j \leq w_n} \prod_{j = 1}^n \frac{1}{\sqrt{x_j}} \bigg [ \int_{\tau = 0}^{\infty} \tau^{-n/2} e^{-[\frac 12 \norm{\theta_0}_2^2 + \sum x_j  ]/(2 \tau)} d\tau \bigg] d \bfx \label{eq:glob_ub_nc1} \\
& = M D^n 2^{n/2-1} \Gamma(n/2-1) w_n^{n/2} \int_{\sum x_j \leq 1} \frac{1}{(\frac 12\norm{\theta_0}_2^2 + w_n \sum x_j)^{n/2-1} } \prod_{j = 1}^n \frac{1}{\sqrt{x_j}}  \, d \bfx \label{eq:glob_ub_nc2} \\ 
& = M D^n 2^{n/2-1} \Gamma(n/2-1) w_n^{n/2} \frac{\Gamma(1/2)^n}{\Gamma(n/2)} \int_{x=0}^1 \frac{x^{n/2 - 1} }{(\frac 12\norm{\theta_0}_2^2 +  w_n x)^{n/2 - 1} } dx.  \label{eq:glob_ub_nc3}
\end{align} 
In the above display, \eqref{eq:glob_ub_nc2} - \eqref{eq:glob_ub_nc3} follows from applying Lemma \ref{lem:dir_formula}  with $h(t) = 1/(\frac 12\norm{\theta_0}_2^2 +  w_n t)^{n/2-1}$. Simplifying constants in \eqref{eq:glob_ub_nc3} as before and using $t/(a+t)$ is an increasing function in $t > 0$ for fixed $a > 0$, we complete the proof by bounding \eqref{eq:glob_ub_nc3} above by
\be
\frac{C w_n}{(n/2-1)} \int_{x=0}^1 \frac{(w_n x)^{n/2 - 1} }{(\frac 12\norm{\theta_0}_2^2 +  w_n x)^{n/2 - 1} }dx &\leq  \frac{C w_n}{(n/2-1)} \bigg(\frac{w_n}{w_n + \frac 12\norm{\theta_0}_2^2} \bigg)^{n/2-1} \\
&\leq \frac{C w_n}{(n/2-1)}  \bigg(\frac{w_n}{\frac 12\norm{\theta_0}_2^2} \bigg)^{n/2-1}. 
\ee
The right hand side of the above display can be bounded above by $e^{- c n \log a_n}$ for some constant $c > 0$. 

For the lower bound on the prior concentration in the non-centered case, we combine \eqref{eq:main_bd1} in the reverse direction along with \eqref{eq:main_bd2}. We then use the same idea as in the centered case to restrict the integral over $\tau$ to $(0, 1)$ in \eqref{eq:glob_lb_nc1}. 
Thus we have
\begin{align}
& \bbP( \norm{\theta - \theta_0}_2 \leq t_n) \geq  K D^n \int_{\sum x_j \leq w_n} \prod_{j = 1}^n \frac{1}{\sqrt{x_j}} \bigg [ \int_{\tau = 0}^{1} \tau^{-n/2} e^{-[\norm{\theta_0}_2^2 + \sum x_j  ]/(2 \tau)} d\tau \bigg] d \bfx. \label{eq:glob_lb_nc1}
\end{align}
Noting that $\norm{\theta_0}_2^2 + \sum x_j \leq \norm{\theta_0}_2^2 + w_n = o(n)$, we can invoke Lemma \ref{lem:incomplete_gamma} to lower bound the inner integral over $\tau$ by $\xi_n \Gamma(n/2-1) 2^{n/2-1}/(\norm{\theta_0}_2^2 + \sum x_j)^{n/2-1} $ and proceed to obtain the same expressions as in \eqref{eq:glob_ub_nc2} \& \eqref{eq:glob_ub_nc3} with $M$ replaced by $K \xi_n$. The proof is then completed by observing that the resulting lower bound can be further bounded below as follows:
\begin{align*} 
& \frac{C w_n}{(n/2-1)} \int_{x=0}^1 \frac{(w_n x)^{n/2 - 1} }{(\norm{\theta_0}_2^2 +  w_n x)^{n/2 - 1} }dx \geq \frac{C w_n}{(n/2-1)} \int_{x=1/2}^1 \frac{(w_n x)^{n/2 - 1} }{(\norm{\theta_0}_2^2 +  w_n x)^{n/2 - 1} }dx  \\
&  \geq \frac{C w_n}{(n/2-1)}  \bigg(\frac{w_n/2 }{(\norm{\theta_0}_2^2 +  w_n/2) }\bigg)^{n/2-1} \geq \frac{C w_n}{(n/2-1)} \bigg(\frac{w_n/2 }{2 \norm{\theta_0}_2^2}\bigg)^{n/2 - 1},
\end{align*}
where the last inequality uses $t_n \leq \norm{\theta_0}_2$ so that $\norm{\theta_0}_2^2 + w_n \leq 2 \norm{\theta_0}_2^2$.
This finishes the proof of Theorem \ref{thm:conc_g}. \qed

\subsection*{ {\bf Proof of Theorem 4.2}}
Let $m_n = (n-3)/2$. Fix $A > 0$.  We set $t_n = s_n$, where $s_n = \sqrt{A \log n}$  so that $w_n = s_n^2 = A \log n$. Also, let $\norm{\theta_0}_2^2 = 2\pi w_n u_n^2$, where $u_n$ is a slowly increasing sequence; we set $u_n = \log (\log n)$ for future references. Finally let $v_n = r_n^2= \sqrt{m_n}$. With these choices, we proceed to show that \eqref{eq:lb_ratio} holds. 

We first simplify  \eqref{eq:lg_ub_stmt} further. The function $x \to x/x(x+a)$ monotonically increases from $0$ to $1$ for any $a > 0$. Thus the integral in  \eqref{eq:lg_ub_stmt} can be further simplified as follows. For any $T_n > 0$, 
\small
\begin{align}
&\int_{\psi_1 = 0}^{\infty} \frac{\psi_1^{m_n}}{ \big\{ \psi_1 + \norm{\theta_0}_2^2/(2\pi w_n) \big \}^{m_n}  } e^{-\psi_1} d \psi_1 \notag \\
& \leq \int_{\psi_1 = 0}^{T_n} \frac{\psi_1^{m_n}}{ \big\{ \psi_1 + \norm{\theta_0}_2^2/(2\pi w_n) \big \}^{m_n}  } e^{-\psi_1} d \psi_1 + \int_{\psi_1 = T_n}^{\infty} e^{-\psi_1} d \psi_1  \leq \bigg( \frac{T_n}{T_n + u_n^2} \bigg)^{m_n} + e^{-T_n}. \label{eq:comb1}
\end{align}
\normalsize
We choose an appropriate $T_n$ which gives us the necessary bound, namely $T_n = u_n \sqrt{m_n} $. Then, using the fact that $(1 - x)^{1/x} \leq e^{-1}$ for all $x \in (0, 1)$, we have
\footnotesize
\begin{align*}
\bigg( \frac{T_n}{T_n + u_n^2} \bigg)^{m_n}  = \bigg( \frac{\sqrt{m_n}}{\sqrt{m_n} + u_n} \bigg)^{m_n}  = \bigg( 1 - \frac{u_n}{\sqrt{m_n} + u_n} \bigg)^{m_n} \leq e^{- m_n u_n/(\sqrt{m_n} + u_n)} \leq e^{- u_n \sqrt{m_n} /2},
\end{align*}
\normalsize
where for the last part used that $e^{-1/x}$ is an increasing function and $\sqrt{m_n} + u_n \leq 2 \sqrt{m_n}$. Thus, using  \eqref{eq:lg_ub_stmt} and substituting $T_n$ in \eqref{eq:comb1} yields, for a global constant $C_1 > 0$, 
\begin{align}\label{eq:lg_simp_ub}
\bbP(\norm{\theta - \theta_0}_2 \leq s_n) \leq \frac{C_1 w_n }{(n/2-1)}  \, \sqrt{\frac{w_n}{ \norm{\theta_0}_2^2}} \,  \, e^{-u_n \sqrt{m_n}/2} .
\end{align}
Next, again using the fact that $x \to x/x(x+a)$ is monotonically increasing, and choosing $c_2 = \infty$, we simplify the lower bound \eqref{eq:lg_lb_stmt}. We have
\begin{align}
& \int_{\psi_1 = c_1 \norm{\theta_0}^2}^{\infty} \frac{\psi_1^{m_n}}{ \big\{ \psi_1 + \norm{\theta_0}_2^2/(\pi v_n) \big \}^{m_n}  } e^{-\psi_1} d \psi_1 
& \geq \bigg(\frac{v_n}{v_n + C}\bigg)^{m_n} e^{- c_1 \norm{\theta_0}_2^2 }, \label{eqn:ttsob}
\end{align}
for some constant $C > 0$. Since $(1 - x)^{1/x} \geq e^{-2}$ for all $x \in (0, 1/2)$ and $e^{-1/x}$ is an increasing function in $x > 0$, we have, 
\begin{align*}
\bigg(\frac{v_n }{v_n + C}\bigg)^{m_n} \geq e^{- \sqrt{m_n}/2}. 
\end{align*}
Hence, the left hand side of \eqref{eqn:ttsob} is bounded below by $e^{-(\sqrt{m_n} + 2c_1 \norm{\theta_0}_2^2)/2}$, resulting in
\begin{align}\label{eq:lg_simp_lb}
\bbP(\norm{\theta - \theta_0}_2 \leq r_n) \geq \frac{C_2 \xi_n v_n }{(n/2-1)}  \, \sqrt{\frac{v_n}{v_n +  \norm{\theta_0}_2^2}} \,  \, e^{-(\sqrt{m_n} + 2c_1 \norm{\theta_0}_2^2)/2}.
\end{align}

Thus, finally, noting that $u_n \to \infty$, \eqref{eq:lb_ratio} follows since
\begin{align*}
\frac{\bbP(|| \theta - \theta_0 ||_2 < s_n )}{\bbP(|| \theta - \theta_0 ||_2 < r_n )}  \times e^{r_n^2}  \leq C' \frac{w_n^{3/2}}{v_n} e^{C (\sqrt{m_n} + \sqrt{n} + \norm{\theta_0}_2^2)}  \, e^{- u_n \sqrt{m_n}  / 2} \to 0,
\end{align*} where $C, C' > 0$ are constants. This finishes the proof of Theorem \ref{thm:lb_lg}. \qed

\subsection*{ {\bf Proof of Theorem \ref{thm:lb_plugin} }}

As before, we assume $\lambda = 1$ without loss of generality, since it can be absorbed in the constant appearing the sequence $\tau_n$ otherwise.    For $t_n < \frac14 \norm{\theta_0}_2 = \frac14 \abs{\theta_{01}}  $,  using \eqref{eq:anderson_UB}  we obtain 
\scriptsize
\begin{eqnarray*}
\bbP(\norm{\theta - \theta_0} < t_n) \leq D^n \tau_n^{-n/2} w_n^{n/2} \int_{\sum x_j \leq 1}  \prod_{j=1}^n \frac{1}{\sqrt{x_j }} \bigg \{\prod_{j=1}^n \int_{\psi_j = 0}^{\infty} \frac{e^{-\psi_j}}{\sqrt{\psi_j}} \exp \bigg(-\frac{w_n x_j + \theta_{01}^2/2}{2\tau_n \psi_j}  \bigg) d\psi_j\bigg\} dx,  \end{eqnarray*} \normalsize
where $w_n = t_n^2$. Using the fact $\int_{0}^{\infty} \frac{1}{\sqrt{x}} \exp \big\{-\big(\frac{a}{x} + x \big)\big\} dx= \sqrt{\pi} e^{-2 \sqrt{a}}$, we obtain
\begin{align}
& \bbP(\norm{\theta - \theta_0} < t_n) \notag \\
& \leq D^n \pi^{n/2} \tau_n^{-n/2} w_n^{n/2}   \int_{\sum x_j \leq 1} \prod_{j=1}^n \frac{1}{\sqrt{x_j} } \exp \bigg\{-2\sqrt\frac{w_n x_j + \theta_{01}^2/2}{2 \tau_n} \bigg\} dx \notag \\
& \leq  \bigg(\frac{D^2 \pi w_n}{\tau_n}\bigg)^{n/2} e^{-\frac{\abs{\theta_{01}}}{\sqrt{\tau_n}}} \int_{\sum x_j \leq 1} \prod_{j=1}^n \frac{1}{\sqrt{x_j} } dx = \bigg(\frac{D^2 \pi w_n}{\tau_n}\bigg)^{n/2} e^{-\frac{\abs{\theta_{01}}}{\sqrt{\tau_n}}} \frac{\Gamma(1/2)^n}{\Gamma(n/2+1)}, \label{eq:plug_in_ub}
\end{align}
where the second to third inequality uses $x_j \geq 0$ and the last integral follows from Lemma \ref{lem:dir_formula}. 
Along the same lines, 
\begin{align}
\bbP(\norm{\theta - \theta_0} < r_n) &\geq \bigg(\frac{D^2 \pi v_n}{\tau_n}\bigg)^{n/2} \int_{\sum x_j \leq 1} \prod_{j=1}^n \frac{1}{\sqrt{x_j} } \exp \bigg\{-2\sqrt\frac{v_n x_j + \theta_{01}^2}{2 \tau_n} \bigg\} dx \notag \\
& \geq \bigg(\frac{D^2 \pi v_n}{\tau_n}\bigg)^{n/2} e^{-\frac{\sqrt{2} (\abs{\theta_{01}} + \sqrt{n v_n})}{\sqrt{\tau_n}}} \frac{\Gamma(1/2)^n}{\Gamma(n/2+1)},  \label{eq:plug_in_lb}
\end{align}
where $v_n = r_n^2/4$. From the second to third equation in the above display, we used $\sqrt{a+b} \leq \sqrt{a} + \sqrt{b}$ and $\sum_{j=1}^n  \sqrt{x_j} \leq \sqrt{n}$ by Cauchy-Schwartz inequality if $x \in \Delta^{n-1}$. Thus, from \eqref{eq:plug_in_ub} \& \eqref{eq:plug_in_lb}, the ratio in \eqref{eq:lb_ratio} can be bounded above as:
\begin{align*}
\frac{\bbP(\norm{\theta - \theta_0} < t_n)}{ \bbP(\norm{\theta - \theta_0} < r_n)}  \leq \bigg(\frac{w_n}{v_n}\bigg)^{n/2} e^{\frac{\sqrt{2 v_n n} + (\sqrt{2} -1) \abs{\theta_{01}}}{\tau_n}}.
\end{align*}
Choose $t_n = s_n$, $r_n = 2 \sqrt{2} s_n$ so that $v_n = 2 w_n = 2 q_n \log(n/q_n)$ and $(w_n/v_n)^{n/2} = e^{- C n}$. Clearly $v_n n /\tau_n \leq C n q_n (\log n)^2$ and hence, $e^{\sqrt{2 v_n n/\tau_n}} = o(e^{C n})$ by assumption. Also, $\abs{\theta_{01}}/\tau_n = o(n)$. Thus, the right hand side of the above display $\to 0$, proving the assertion of the Theorem.

\section*{Acknowledgements}
NSP is indebted to Prof. Gine for generously sharing his time and ideas.
We thank Professors Sara van de Geer and Richard Nickl for encouragement and the referees for constructive criticism.
  \bibliographystyle{elsarticle-num} 
 \bibliography{cov_refs}


%
%
%
\end{document}